\documentclass[11pt]{amsart}

\usepackage{amsmath,amsfonts,amssymb}
\usepackage{array,delarray}
\usepackage[dvips]{graphicx}
\usepackage{dsfont}
\usepackage[square]{natbib}

\newcommand{\reviewtimetoday}[2]{
}
\reviewtimetoday{\today}{Final Version}

\newcommand{\zz}{\mathds{Z}}

\newcommand{\ff}{\mathds{F}}

\newcommand{\supp}{\mathrm{supp}}

\newtheorem{theorem}{Theorem}[section]   % Numbered within each section
\newtheorem{corollary}[theorem]{Corollary}     % Numbered along with thm
\newtheorem{lemma}[theorem]{Lemma}         % Numbered along with thm
\newtheorem{proposition}[theorem]{Proposition}  % Numbered along with thm

\theoremstyle{definition}
\newtheorem{definition}[theorem]{Definition}   % Numbered along with thm

\theoremstyle{remark}
\newtheorem{example}[theorem]{Example}        % Numbered along with thm

\numberwithin{equation}{section}     % Number equations within sections

\begin{document}

\title[Monomial Dynamical Systems over Finite Fields]{Monomial Dynamical Systems over Finite Fields}
\thanks{The second and third authors were supported in part by NSF Grant Nr. DMS-0511441.
The third author was also partially supported by NIH Grant Nr. RO1 GM068947-04.}
\author{O. Col\'on-Reyes}
\address{Mathematics Department,
University of Puerto Rico at Mayag\"{u}ez,
Mayag\"{u}ez, PR 00681}
\email{ocolon@math.uprm.edu}
\author{A. S. Jarrah}
\address{Virginia Bioinformatics Institute,
Virginia Tech,
Blacksburg, VA 24061-0477, USA}
\email{ajarrah@vbi.vt.edu}

\author{R. Laubenbacher}
\address{Virginia Bioinformatics Institute,
Virginia Tech,
Blacksburg, VA 24061-0477, USA}
\email{reinhard@vbi.vt.edu}

\author{B. Sturmfels}
\address{Department of Mathematics, University of California at Berkeley,
Berkeley, CA 94720}
\email{bernd@math.berkeley.edu}


\begin{abstract}
An important problem in the theory of finite dynamical systems is to
link the structure of a system with its dynamics.  This paper
contains such a link for a family of nonlinear systems over an
arbitrary finite field. For systems that can be described by
monomials, one can obtain information about the limit cycle
structure from the structure of the monomials. In particular, the
paper contains a sufficient condition for a monomial system to have
only fixed points as limit cycles.  The condition is derived by
reducing the problem to the study of a Boolean monomial system and a
linear system over a finite ring.
\end{abstract}

\maketitle

\section{Introduction}
%=====================
 Finite dynamical systems are time-discrete dynamical
systems on finite state sets. Well-known examples include cellular
automata and Boolean networks, which have found broad applications
in engineering, computer science, and, more recently,
computational biology.  (See, e.g., \cite{K, AO, CS1, LS} for
biological applications.) More general multi-state systems have
been used in control theory \cite{GGP, LBL, ML1, ML2}, the design
and analysis of computer simulations \cite{BR, BMR1, BMR2, LP2}.
One underlying mathematical question that is common to many of
these applications is how to analyze the dynamics of the models
without actually enumerating all state transitions, since
enumeration has exponential complexity in the number of model
variables.  The present paper is a contribution toward an answer
to this question.

For our purposes, a finite dynamical system is a function
$f:X \longrightarrow X$, where $X$ is a finite set \cite{LP1}.
The dynamics of $f$ is generated by iteration of $f$ and is
encoded in its \emph{phase space} $\mathcal S(f)$,
which is a directed graph defined as follows.
The vertices of $\mathcal S(f)$ are the elements of
$X$. There is a directed edge $a\rightarrow b$ in $\mathcal S(f)$
if $f(a)=b$.  In particular, a directed edge from a vertex to itself
is admissible. That is, $\mathcal S(f)$ encodes all state
transitions of $f$, and has the property that every vertex has
out-degree exactly equal to $1$.  Each connected graph component of
$\mathcal S(f)$ consists of a directed cycle, a so-called {\em limit
cycle}, with a directed tree attached to each vertex in the cycle,
consisting of the so-called {\em transients}.

Any Boolean network can be viewed as a  finite dynamical system
$f: \ff_2^n \longrightarrow \ff_2^n$, where $\ff_2$ is the finite
field on two elements and $n \geq 1$. In this paper, we study finite
dynamical systems $f : \ff_q^n \longrightarrow \ff_q^n$, where $\ff_q$ is any finite
field and $|\ff_q| = q$.
To be precise, we present a family of
nonlinear finite systems for which the above question can be answered, that
is, for which one can obtain information about the dynamics from the
structure of the function.

Let  $f:\ff_q^n \longrightarrow  \ff_q^n$, $n\geq 1$ be a finite
dynamical system. Observe that $f$ can be described in terms of
its coordinate functions $f_i:\ff_q^n \longrightarrow  \ff_q$,
that is, $f=(f_1,\ldots ,f_n)$.  It is well known that any
set-theoretic function $g:\ff_q^n \longrightarrow  \ff_q$ can be
represented by a polynomial in $\ff_q[x_1,\ldots ,x_n]$, see
\cite[pp. 369]{LN}.  This polynomial can be chosen uniquely so
that any variable in it appears to a degree less than $q$. That
is, for any $g$ there is a unique $\bar h\in \ff_q[x_1,\ldots
,x_n]/\langle x_i^q-x_i|i=1,\ldots ,n\rangle$, such that
$g(a)=h(a)$ for all $a\in \ff_q^n$. Consequently, any finite
dynamical system over a finite field can be represented as a
polynomial system.  This is the point of view we take in this
paper.

In the case that all the $f_i$ are linear polynomials without constant
term,  the dynamics of the \emph{linear system} $f$  can completely be determined
from its matrix representation \cite{He}. Let $A$ be a
matrix representation of a linear system $f:\ff_q^n \longrightarrow \ff_q^n$.
Then the number of limit cycles and their length, as well as
the structure of the transients, can be determined from the
factorization of the characteristic polynomial of the matrix $A$.
The structure of the limit cycles had been determined earlier by
Elspas \cite{El}, and for affine systems by Milligan and Wilson
\cite{MW}.

In this paper we focus on the class of nonlinear systems described
by special types of polynomials, namely \emph{monomials}. That is,
we consider systems $f=(f_i)$, so that each $f_i$ is a polynomial
of the form $x_{i_1}^{a_{i1}}x_{i2}^{a_{i2}}\cdots
x_{i_r}^{a_{ir}}$, or a constant.  Without loss of generality we
can assume that no coordinate function is constant, since the
general case is easily reduced to this. Some classes of monomial
systems and their dynamic behavior have been studied before:
Monomial cellular automata \cite{BG,Kari}, Boolean monomial
systems \cite{CLP}, monomial systems over the p-adic numbers
\cite{KhrennikovNilsson:01,Nilsson:03}, and monomial systems over
finite fields \cite{Omar,VasigaShallit}.

%Associated to a general polynomial
%system one can construct its {\em dependency graph} $\mathcal D(f)$,
%whose vertices $v_1,\ldots ,v_n$ correspond to the variables of the
%$f_i$.  There is a directed arrow $v_i\rightarrow v_j$ if $x_j$
%appears in $f_i$.  For Boolean monomial systems the dependency graph
%in fact allows the unambiguous reconstruction of the system.

In \cite{CLP} a special class of Boolean monomial systems was
studied, namely those which have only fixed points as limit cycles,
so-called \emph{fixed point systems}.  The motivation for considering this
class is the use of polynomial systems as models for biochemical
networks. Depending on the experimental system considered, such
networks often exhibit steady state dynamics. That is, their dynamic
models have phase spaces whose limit cycles are fixed points.  For
the purpose of model selection it would be useful to have a
structural criterion to recognize fixed point systems.  The main
result of the present paper is to reduce this question for monomial
systems over a general finite field $\ff_q$ to the
same question for an associated Boolean monomial system and a linear
system over a ring of the form $\zz/(q-1)$.
%One can then apply the
%criterion in \cite{CLP} to the Boolean system.
%%%%%%%%%%%%%%%%%%%%%%%%%%%%%%%%%%%%%%%%%%%%%%%%%%%
%We derive an analogous necessary condition for linear
%systems over the more general rings $\mathbf Z/(q-1)$.

\section{Reduction of monomial systems}
%=======================================
Let
$f=(f_1,\ldots ,f_n):\ff_q^n\longrightarrow \ff_q^n$ be a polynomial system,
with each $f_i$ a monomial, that is,
$f_i=x_1^{a_{i1}}x_2^{a_{i2}}\cdots x_n^{a_{in}}$, with $a_{ij}$ a
nonnegative integer.  That is, $f$ can be described by the exponent
matrix $A=(a_{ij})$.  We first associate with $f$ a Boolean monomial
system $T(f)$ and a linear system $L(f)$ over the ring $R=\zz/(q-1)$.
 Recall from \cite{CLP} that $f$ is called a
{\it fixed point system} if all limit cycles of $f$ consist of a
fixed point. We will show that $f$ is a fixed point system if and
only if $T(f)$ and $L(f)$ are fixed point systems.

\begin{definition}
For $\textbf{u} = (u_1,\dots,u_n) \in \ff_q^n$, we define support
$\textbf{u}$, denoted by $\supp(\textbf{u})$, to be $\textbf{v} = (v_1,\dots,v_n)$, where
\[
v_i = \left\{
\begin{array}{cc}
  1 & \mbox{ if  $u_i \neq 0$}; \\
  0 & \mbox{ if  $u_i = 0$}. \\
\end{array}
\right.
\]
\end{definition}

The monomial system $f =(f_1,\dots,f_n)$ induces a Boolean
monomial system $T(f) = (g_1,\dots,g_n)$ on $\ff_2^n$ by setting
$g_i= x_1^{v_{i1}}\cdots x_n^{v_{in}}$, where $f_i =
x_1^{u_{i1}}\cdots x_n^{u_{in}}$ and $\textbf{v} =
\supp(\textbf{u})$.

\begin{lemma}\label{T}
There is a commutative diagram
$$
\begin{matrix}
\ff_q^n         & \stackrel{f}{\longrightarrow} &\ff_q^n\\
\downarrow  &                   &\downarrow\\
\ff_2^n& \stackrel{T(f)}{\longrightarrow}   &\ff_2^n
\end{matrix}
$$
\end{lemma}

\begin{proof}
The proof is a straightforward verification, since
$$
\supp(f(\mathbf u))= f(\supp(\mathbf u)).
$$
\end{proof}
Since $f = T(f)$ on the set of all $\textbf{u}$ such
that $\supp(\textbf{u}) = \textbf{u}$, the following
corollaries are straightforward.
\begin{corollary}\label{cor-1}
The phase space of $T(f)$ is a subgraph of the phase space of $f$.
\end{corollary}
\begin{corollary}
Suppose that $T(f)$ is a fixed-point system.
If $\{\mathbf u, f(\mathbf u), \dots, f^t(\mathbf u) = \mathbf u \}$
is a cycle in the phase space of $f$, then
$\supp(\mathbf u) = \supp(f^i(\mathbf u))$ for all $ 1 \leq i \leq t$.
\end{corollary}
For more results in this direction, see \cite{Omar}.

\begin{figure}[htbp]
%\label{fig:1}
\begin{center}
\includegraphics[width=6in]{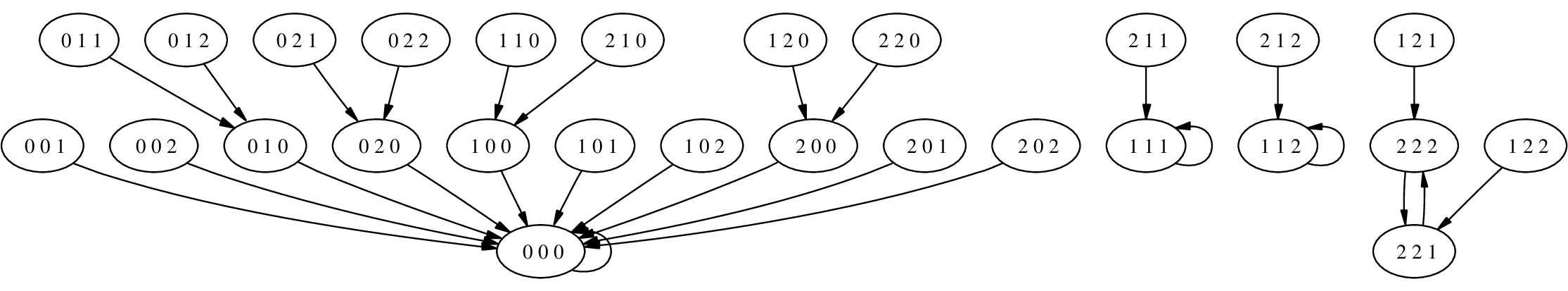}
\caption{The phase space of $f$.}
\label{fig:1}
\end{center}
\end{figure}

\begin{figure}[htbp]
\label{fig:2}
\begin{center}
\includegraphics[width=1.8in]{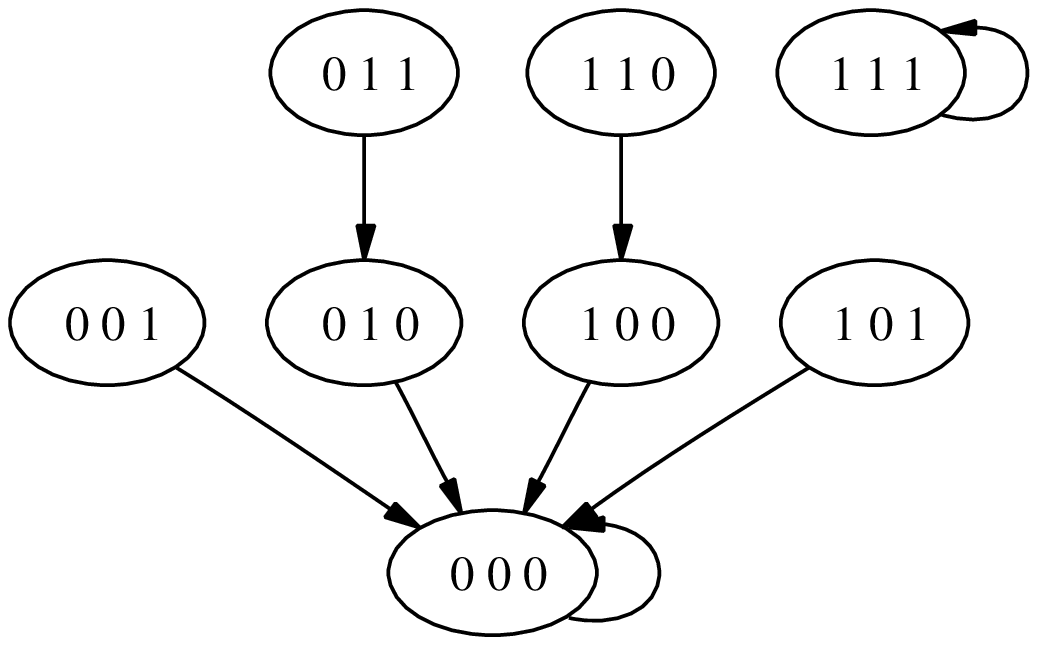}
\caption{The phase space of $T(f)$.}
\end{center}
\end{figure}

\begin{example}\label{ex-1}
Let
$$
f=(x_1^2x_2, x_2x_3^2, x_1^2x_2x_3):\ff_3^3\longrightarrow \ff_3^3.
$$
\end{example}
Figures \ref{fig:1} and \ref{fig:2} display the phase spaces of the
system $f$  and its "Booleanization" $T(f)$, respectively.
%It is
%easy to see the projection of phase spaces that is induced by the
%support map projection of $\ff_3^3$ onto $\ff_2^3$.

Next we associate to $f$ an $n$-dimensional linear system over a
finite ring.  Observe first that $\ff_q^*=\ff_q-\{0\}$ is isomorphic, as an
Abelian group, to $\zz/(q-1)$ via an isomorphism
$$
log: \ff_q^* \longrightarrow \zz/(q-1),
$$
given by the choice of a generator for the cyclic group $\ff_q^*$.
Note first that the set of vectors $\mathbf u\in \ff_q^n$ with all
non-zero entries is invariant under $f$.

Let $\alpha$ be a generator for the cyclic group $\ff_q^*\cong \mathbf
\zz/(q-1)$, and let
$$
\mathbf u=(u_1,\ldots ,u_n)=(\alpha^{s_1},\ldots
,\alpha^{s_n})=\alpha^{\mathbf s}.
$$
Then $f(\mathbf u)=f(\alpha^{\mathbf s})=\alpha^{A\cdot \mathbf s}$.

\begin{definition}
Define $L(f):(\zz/(q-1))^n\longrightarrow (\zz/(q-1))^n$ by
$$
L(f)(\alpha^{\mathbf s})=\alpha^{A\cdot \mathbf s}.
$$

\end{definition}

As defined, $L(f)$ is a linear transformation of $\zz$-modules.
But we can consider it as a linear transformation of
$\zz/(q-1)$-modules, viewing $\zz/(q-1)$ as a (finite) ring, which
we denote by $R$.  That is, we have the linear transformation
$$
L(f)=A:R^n\longrightarrow R^n.
$$

The proof of the following lemma is a straightforward
verification.

\begin{lemma}\label{L}
There is a commutative diagram
$$
\begin{matrix}
(\ff_q^*)^n& \stackrel{f}\longrightarrow &(\ff_q^*)^n\\
\downarrow&&\downarrow\\
R^n& \stackrel{L(f)}\longrightarrow &R^n
\end{matrix}
$$
\end{lemma}
Note that the vertical arrows are isomorphisms.
This implies that they preserve the phase space structure,
including the length of limit cycles \cite{LP1}. In particular,
we have the following corollary.
\begin{corollary}\label{cor-2}
The phase space of $L(f)$ is isomorphic to the subgraph of the
phase space of $f$ consisting of all states with support vector
$(1,1,\dots,1)$.
\end{corollary}

\medskip

\begin{figure}[htbp]
%\label{fig:3}
\begin{center}
\includegraphics[width=2in]{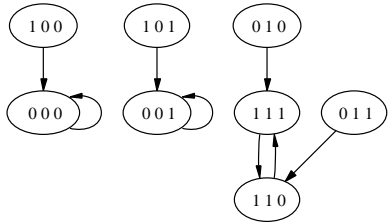}
\caption{The phase space of $L(f)$.}
\label{fig:3}
\end{center}
\end{figure}

\noindent \emph{Example }\ref{ex-1} (Cont.). For the monomial system
$f$ in Example \ref{ex-1}, $L(f) :(\zz/2)^3\longrightarrow (\zz/2)^3$
is defined by $L(f)(\mathbf s) = A\cdot \mathbf s$, where
\[
A =
\left(%
\begin{array}{ccc}
  0 & 1 & 0 \\
  0 & 1 & 0 \\
  0 & 1 & 1 \\
\end{array}%
\right).
\]
The phase space of $L(f)$ is given in Figure \ref{fig:3}.

We now come to the main result of this section.
\begin{theorem}
Let $f:\ff_q^n\longrightarrow \ff_q^n$ be a monomial dynamical system.
Then $f$ is a fixed point system if and only if $T(f)$ and $L(f)$
are fixed point systems.
\end{theorem}

\begin{proof} From Corollaries \ref{cor-1} and \ref{cor-2}, if $f$ is a fixed
point system, then so are $T(f)$ and $L(f)$. For the converse, we
assume that $L(f)$ and $T(f)$ are fixed point systems, but $f$ is
not. For each limit cycle of $f$, either all states involved have
all coordinates non-zero, or all states involved have at least one
zero coordinate.  In the first case it follows that $L(f)$ has a
limit cycle of the same length.  Hence, if $f$ has a limit cycle
of length greater than $1$, then it can involve only states with
at least one zero coordinate.

Let $\mathbf a_1,\ldots ,\mathbf a_t$ be the states in the limit
cycle.  Since this limit cycle has to project to a fixed point of
$T(f)$ it follows that the $\mathbf a_i$ have the same support
vector, i.e., the same pattern of zero entries, and differ only in
the non-zero coordinates. Furthermore, the monomials in the
nonzero coordinates do not involve any variables corresponding to
zero coordinates. Thus, if we construct new states $\mathbf b_i$
by replacing each $0$ in $\mathbf a_i$ by a $1$, the $\mathbf b_i$
will be part of a limit cycle of length at least $t$, which is a
contradiction. This completes the proof of the theorem.
\end{proof}

\section{Linear Systems over Finite Commutative Rings}
%====================================================
The theorem in the previous section shows that to decide whether a
given monomial system $f$, over a finite field $\ff_q$,
 is a fixed point system it is sufficient to
decide this question for an associated Boolean system, for which a
criterion has been developed in \cite{CLP}, and a certain linear
system over a finite ring $\zz/m$.  It therefore remains to
develop a criterion for linear systems over finite commutative
rings that helps decide whether the system is a fixed point
system.  Here we reduce the case of a general $\zz/m$ to that of
$m$ being a prime power.

Let $m = st$ for relatively prime integers $s$ and $t$, and let $f$
be a linear system over $\zz/m$ of dimension $n$.   Choosing
an isomorphism $\zz/m \cong \zz/s\times \zz/t$ we
see that $f$ is isomorphic to a product
$$
g \times h :(\zz/s)^n\times (\zz/t)^n \longrightarrow (\zz/s)^n\times (\zz/t)^n,
$$
where $g$ and $h$ are linear systems over $\zz/s$ and $\zz/t$, respectively.
Using the fact that the phase space of $f$ is then the direct
product, as directed graphs, of the phase spaces of $g$ and
$h$ (see, e.g., \cite{He} or \cite{JLV}), we obtain the following
result.

\begin{proposition}
Let $m =st$ for relatively prime integers $s$ and $t$, and let $f$
be a linear system over $\zz/m$ of dimension $n$. Let $g$
and $h$ be the induced linear transformations over $\zz/s$
and $\zz/t$, respectively.  Then $f$ is a fixed point system
if and only if $g$ and $h$ are fixed point systems.
\end{proposition}

For the purpose of developing a criterion to recognize fixed point
systems it is therefore sufficient to study linear systems over
rings of the form $\zz/p^r$ for primes $p$.  The following theorem
provides a criterion for a further reduction of the problem to a
linear system over the prime field $\zz/p$.

\begin{figure}[htbp]
\begin{center}
\includegraphics[width=6in]{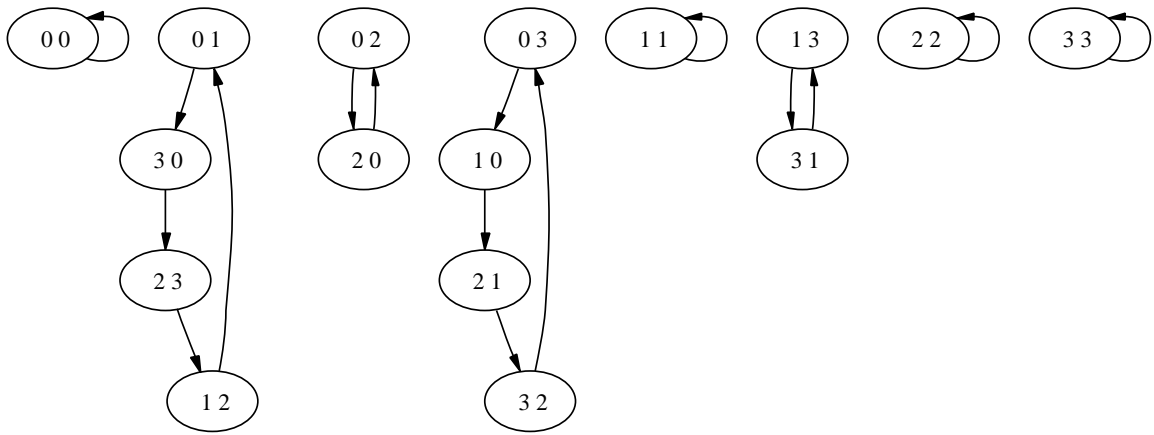}
\caption{The phase space of $f$.}
\label{fig:4}
\end{center}
\end{figure}

\begin{figure}[htbp]
\begin{center}
\includegraphics[width=1.5in]{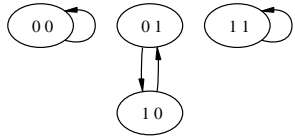}
\caption{The phase space of $g$.}
\label{fig:5}
\end{center}
\end{figure}

\begin{theorem}
Let $f = (f_1,\dots,f_n) : (\zz/p^r)^n \longrightarrow
(\zz/p^r)^n$ be a linear map, and let $g$ be the projection map
of $f$ on $\zz/p$. That is
$g= (g_1,\dots,g_n): (\zz/p)^n \longrightarrow (\zz/p)^n$,
where $g_i = f_i \mod p$.
Then the phase space of $g$ is isomorphic to a subgraph of the
phase space of $f$.
\end{theorem}
\begin{proof}
Let $\sigma : (\zz/p)^n \longrightarrow (\zz/p^r)^n$
be given by $\sigma(\textbf{a})= \sigma(a_1,\dots,a_n) =
(a_1p^{r-1}, \dots, a_n p^{r-1})= \textbf{a}p^{r-1}$. Then it is
easy to check that $\sigma \circ g = f \circ \sigma$, since the
$f_i$ are linear maps for all $i$. Therefore, it is
straightforward to check that $g(\textbf{a}) = \textbf{b}$ if and
only if $f(\textbf{a}p^{r-1}) = \textbf{b}p^{r-1}$, and hence the
phase space of $g$ is isomorphic to subgraph of the phase space of
$f$.
\end{proof}
\begin{corollary}
Let $f$ and $g$ be as above. If $g$ is not a fixed point system,
then $f$ is not a fixed point system.
\end{corollary}
The dynamics of projection maps has been studied in \cite{Omar}.
\begin{example}
Let $f: (\zz/4)^2\longrightarrow (\zz/4)^2$ be given
by $f = (2x_1+3x_2, x_1)$.
Then $g = (x_2,x_1) :(\zz/2)^2\longrightarrow (\zz/2)^2$.
The phase spaces of $f$
and $g$ are in Figures \ref{fig:4} and \ref{fig:5}, respectively.
\end{example}

It remains to study the dynamics of linear systems over finite
rings, in particular to find a criterion for a linear system to be
a fixed point system.  Together with the results in this paper,
such a criterion would provide an algorithm for deciding whether a
monomial system over an arbitrary finite field is a fixed point
system.  However, it appears to be a difficult problem to
understand the dynamics of linear systems even over rings of the
form $\zz/p^r$, due to the lack of unique factorization in the
polynomial ring $\zz/p^r[x]$.  See, e.g., \cite{B}.

%%%%%%%%%%%%%%%%%%%%%%%%%%%%%%%%%%%%%%%%%%%%%%%%%%%%%%%%%%%%%%
%\bibliographystyle{aiaa}
%\bibliographystyle{klunum}
%\bibliographystyle{apa}
%\bibliographystyle{siam}
%\bibliography{../../ext_bib_2}

%%%%%%%%%%%%%%%%%%%%%%%%%%%%%%%%%%%%%%%%%%%%%%%%%%%%%%%%%%%%%%
\end{document}